%% file: agt-5-20.tex
\documentclass{gtart_h}
\input agtout

\lognumber{20}
\volumenumber{5}
\volumeyear{2005}
\papernumber{20}
\pagenumbers{443}{462}
\received{14 September 2004} 
\revised{29 March 2005}
\accepted{8 April 2005}
\published{26 May 2005}

\usepackage{amssymb, graphicx,psfrag}
\def\psfraga <#1,#2>#3#4{%
\psfrag {#3}{\smash{\rlap{\kern #1 \raise #2\hbox{#4}}}}}

\def\tl{\triangleleft}
\def\tr{\triangleright}

\newtheorem{theorem}{Theorem}[section]
\newtheorem{proposition}[theorem]{Proposition}

\theoremstyle{definition}
\newtheorem{definition}{Definition}[section]
\newtheorem*{remark}{Remark}

\author{Sam Nelson} 
\address{University of California, Riverside\\900 University Avenue, 
Riverside, CA 92521, USA}
\email{knots@esotericka.org}
\title{Signed ordered knotlike quandle presentations}


\begin{abstract}
We define enhanced presentations of quandles via generators and
relations with additional information comprising signed operations and
an order structure on the set of generators.  Such a presentation
determines a virtual link diagram up to virtual moves.  We list formal
Reidemeister moves in which Tietze moves on the presented quandle are
accompanied by corresponding changes to the order structure.  Omitting
the order structure corresponds to replacing virtual isotopy by welded
isotopy.
\end{abstract}

\begin{document}

\maketitle

{\small\it Dedicated to the memory of my mother, Linda M Nelson, who
was killed in an automobile accident while this paper was in
preparation.\leftskip25pt\rightskip25pt\par}

\section{Introduction}

The knot quandle was defined in \cite{J} and has been studied in a
number of subsequent works. It is well known that the knot quandle is
a classifying invariant for classical knots, in the sense that the
quandle contains the same information as the fundamental group system
and hence determines knots and links up to (not necessarily
orientation-preserving) homeomorphism of pairs. Whether the quandle
alone is a complete invariant of knot type, however, depends on the
meaning of ``equivalence'' -- in particular, if ``equivalent'' means
``ambient isotopic'', then the quandle alone does not classify
knots. Indeed, there are examples of pairs of classical (and virtual)
knots with isomorphic knot quandles which are not ambient (or
virtually) isotopic, such as the left and right hand trefoil knots.
This shows that a quandle isomorphism need not translate to a
Reidemeister move sequence. Moreover, Reidemeister moves on a knot
diagram correspond to sequences of Tietze moves on the presented
quandle, but the converse is not generally true. We wish, then, to
characterize which Tietze moves on quandle presentations do correspond
to Reidemeister moves and hence preserve ambient (or virtual) isotopy
class. We will assume that all knots are endowed with a choice of
orientation, ie, we will concern ourselves only with oriented knots
and links.

Examples of nontrivial virtual knots with trivial knot quandle are
given in \cite{D}.  Virtual knots are equivalence classes, under the
equivalence relation generated by the three Reidemeister moves,
of 4--valent graphs (both planar and non-planar) with vertices interpreted as 
crossings; classical knot theory may be regarded as the special case of 
virtual knot theory in which we restrict our attention to planar graphs. 
Virtual knots have two associated quandles, an \textit{upper quandle} defined 
in the usual combinatorial way (described in section 3) while ignoring any 
virtual crossings, and a \textit{lower quandle}, defined as the upper quandle 
of the knot diagram obtained by ``flipping over'' the original knot diagram 
by taking a mirror image and switching all under/overcrossings. If the knot 
is classical, the ``flipping over'' operation is an ambient isotopy, and the 
resulting upper and lower quandles are isomorphic, though the quandle 
presentations defined by the knot diagrams may bear little obvious 
resemblance. If the knot is not classical, however, the upper and lower 
quandles are typically distinct; see \cite{GPV}.

One approach to the problem of finding a complete algebraic invariant
for virtual knots involves combining the upper and lower quandles into
a single algebraic structure, namely the biquandle (defined in
\cite{FRS:species} and applied to virtual knots in \cite{FJK}).  In this
paper we take a different approach, defining additional structure on
the usual presentation of the upper quandle which permits
reconstruction of a virtual knot diagram from the enhanced
presentation.

In section \ref{sec:qp}, we recall presentations of quandles by
generators and relations, noting that every finite quandle has a
presentation resembling that determined by a knot or link diagram. In
section \ref{sec:kp} we recall some facts about virtual knots and the
definition of the knot quandle, then present our primary definition,
signed ordered knotlike quandle presentations. In section
\ref{sec:frm}, we examine how the Reidemeister moves affect the order
structure, defining formal Reidemeister moves on SOKQ
presentations.  In section \ref{sec:weld}, we show how to apply our
results to welded isotopy classes (also known as weak virtual isotopy
classes) and in section \ref{sec:framed} extend to cover the framed
case.

\section{Quandle presentations}\label{sec:qp}

\newcounter{ax}
A \textit{quandle} is a set $Q$ with a non-associative binary operation 
$\tr\co  Q\times Q\to Q$ satisfying 
\begin{list}{(q\roman{ax})}{\usecounter{ax}}
\item for every $a\in Q$, we have  $a\tr a=a$,
\item for every $a,b\in Q$ there exists a unique $c\in Q$ with
$a=c\tr b$, and
\item for every $a,b,c \in Q$, we have 
$(a\tr b)\tr c=(a\tr c)\tr(b\tr c)$.
\end{list}
As expected, a \textit{homomorphism} of quandles is map $f\co Q\to Q'$
such that \[f(x\tr y) = f(x)\tr f(y),\] and a bijective homomorphism
of quandles is an \textit{isomorphism} of quandles.  The quandle
operation is asymmetrical; it is an action of the set $Q$ on itself.
Several authors have written on quandles, and a number of different
notational styles are in common use:  \cite{FR} uses exponential
notation, where $a\tr b$ is denoted $a^b$, while \cite{CJKLS} uses an
asterisk $*$ in place of the triangle $\tr$. Moreover, some authors
put the action on the right, while others put it on the left.

For the purpose of defining quandles via presentations, it is convenient
to follow \cite{J} using the triangle notation with the action on the right, 
so that $x\tr y$ means ``the result of the action of $y$ on $x$.'' The 
existence and uniqueness requirements of axiom (qii) imply that each quandle 
comes with a second operation $\tl$ satisfying $(x \tl y) \tr y =x$. 
Specifically, axiom (qii) is equivalent to the statement that for all 
$y\in X$, the map $f_y(x)=x\tr y$ is a bijection. We may then denote 
$f^{-1}_y(x)=x\tl y$, and we have $(x\tr y) \tl y = f_y^{-1}(f_y(x)) = x$ and 
$(x\tl y) \tr y = f_y(f_y^{-1}(x)) = x$. This operation $\tl\co  Q\times Q\to Q$ 
itself defines a quandle structure on $Q$, called the \textit{dual} of $Q$. 
A quandle is \textit{self-dual} if it is isomorphic to its dual, that is, if 
there is a bijection $f\co Q\to Q$ with $f(x\tr y) = f(x)\tl f(y)$. The dual
operation $x\tl y$ is also denoted $x^{\bar{y}}$, $x\bar{*} y$, or
$x\tr^{-1} y$.

We will use the symbol $\diamond\in \{\tr, \tl\}$ as a generic quandle 
operation, so the notation $x\diamond y$ can mean either $x\tr y$ or 
$x\tl y$. When we have already used $\diamond$ in a formula, we will use 
$\bar{\diamond}$ to specify the opposite operation, and when we need
to specify several possibly different quandle operations, we will use 
subscripts: $\diamond_1, \diamond_2$ etc.

Although the quandle operations are not associative, we can write any
element of the quandle in the form
\[(\dots (((x_1 \diamond_1 x_2) \diamond_2 x_3 ) 
\diamond_3 \dots)\diamond_{n-1} x_n),\]
$$x\diamond_1 (y\diamond_2 z) = ((x\bar{\diamond}_2 z)\diamond_2
z)\diamond_1 (y\diamond_2 z) = (((x\bar{\diamond}_2 z) \diamond_1 y)
\diamond_2 z).\leqno{\rm since}$$ 
This form is called \textit{fully left assocated} and we shall often
assume, without comment, that words are presented in this form.  The
\textit{length} of such a word is the number of operations ($\tr$ or
$\tl$) which occur.

We refer to \cite{FR} for basic results about presentations of
quandles, free quandles and the Tietze theorem.  \cite{FR} is written
in terms of the more general concept of racks (axiom (qi) is omitted)
but the treatment specialises to quandles in the obvious way.  A
quandle presentation can contain both primary and secondary (or operator)
generators and relations.  Here we shall only be concerned with primary
generators and relations.  Such a presentation has the form $\langle
X\mid R\rangle$, where $X=\{x_1,\ldots,x_n\}$ is the generating set
and $R=\{r_i\sim r_i', i=1\dots m\}$ is the set of relations.  Here
each $r_i$ and $r_i'$ is a word in $x_1,\ldots,x_n$.

\begin{definition}
A relation of the form $x\sim y \diamond z$ where $x, y,$ and $z$ are 
generators is a \textit{short} relation. In a relation of this form, the 
generator $x$ is the \textit{output operand}, $y$ is the \textit{input 
operand}, and $z$ is the \textit{operator}.
\end{definition}

\begin{proposition} \label{sform}
Every finitely presented quandle $\langle X \ | \ R\rangle $ has a 
presentation such that
\begin{list}{$\bullet$}{}
\item every relation is short,
\item at most one of $\left\{ x\sim y\diamond z, 
x'\sim y\diamond z \right\}$ appears in $R$,  
\item at most one of $\left\{ x\sim y\diamond z, 
x\sim y'\diamond z \right\}$ appears in $R$,
\item at most one of $\left\{ x\sim y\tr z, 
y\sim x\tl z \right\}$ appears in $R$,
\item at most one of $\left\{ x\sim y\tr y, 
x\sim y\tl y \right\}$ appears in $R$.
\end{list}
\end{proposition}

\begin{proof}
Every relation has the form $w_1\sim w_2$ where $w_i$ are words in
the generators.  We may assume $w_1,\ w_2$ are not both generators,
since if they are, we can delete one and replace every occurrence of
the deleted generator with the equivalent generator without changing
the quandle.  Each word $w_i$ can be written in fully left associated
form, so we may further assume every relation has the form
\[
x_1 \diamond_1 x_2 \diamond_2 \dots \diamond_{n-1} x_n \ \sim \
y_1 \diamond_{n+1} y_2 \diamond_{n+2} \dots \diamond_{n+m-1} y_m.
\]
If we already have a relation $z\sim x_1\diamond_1 x_2$, we can replace 
$x_1\diamond_1 x_2$ with $z$, shortening the left hand side of the relation; 
if not, we can introduce a new generator $z$ with defining relation 
$z\sim x_1\diamond_1 x_2$ and do the replacement. Repeating the procedure, 
we can reduce the left side (say) to a single generator and the right side 
to a word $u\diamond v$ of length one, so that all relations have
the short form $x\sim y\diamond z$.

If both $x\sim y\diamond z$ and $x'\sim y\diamond z$ are in $R$, then 
$x\sim x'$ and replacing all instances of $x'$ with $x$ and deleting the 
relation $x'\sim y\diamond z$ and the generator $x'$ yields an isomorphic 
quandle. Similarly, since $x\sim y\diamond z$ is equivalent to 
$y\sim x\bar{\diamond} z$, the third claim reduces to the second.

By the quandle laws $x\sim y\tr z$ implies $x\tl z \sim (y\tr z) \tl z$ which
implies $x\tl z \sim y$.  Hence   
the relations $x\sim y\tr z$ and $y\sim x\tl z$ are  equivalent, so 
we can choose one and delete the other without changing the presented quandle.
Similarly, $x\sim y\tr y$ is equivalent to $x\sim y\tl y$.
\end{proof}

\begin{remark}It is also true that $x\sim y\tr y$ is equivalent to $x\sim y$
by (qi) and therefore relations of the form $x\sim y\tr y$ can be
removed from the presentation.  However it is convenient to keep these
redundant relations as they occur naturally in presentations coming
from diagrams.  They are also necessary for the analogous framed case,
where racks are used rather than quandles and axiom (qi) is not
assumed to hold (see section \ref{sec:framed}).  Note also that the
relations $x\sim y\diamond z$ and $x\sim y\diamond z'$ are
independent, in general.\end{remark}

\begin{definition}
A quandle presentation $Q=\langle X \ | \ R \rangle $ of the type described in 
the preceding proposition is in \textit{short from}. 
\end{definition}

A virtual knot or link diagram determines a quandle presentation in short 
form, as we shall see in the next section. In particular, the quandle 
presentations determined by the diagrams before and after a Reidemeister move 
are both in short form, so a sequence of Tietze moves on a knot quandle must 
determine a sequence of short-form presentations in order to correspond to a 
sequence of Reidemeister moves. 

Since every quandle has a short form presentation, it is natural to ask when 
a virtual knot diagram can be reconstructed from a short form quandle 
presentation and to what degree the isomorphism class of the quandle 
determines the resulting virtual knot. We will see that a short form quandle 
presentation meeting certain sufficient conditions, together with 
some additional structure, determines a virtual link diagram.

\section{Knotlike presentations and link diagrams}\label{sec:kp}

Interest in quandles is motivated primarily by their utility in
defining invariants of knots and links, disjoint unions of embedded
circles in $S^3$ or another 3--manifold. A knot or link with a choice
of orientation is an \textit{oriented} knot or link. Oriented knots
and links in $S^3$ may be studied combinatorially as equivalence
classes of \textit{link diagrams}, 4--valent graphs with two edges
oriented in and two oriented out at every vertex, with vertices
interpreted as crossings and enhanced with crossing information, under
the equivalence relation generated by the three familiar Reidemeister
moves. If the underlying graph of a
link diagram is a plane graph, ie, a planar graph actually embedded in
the plane, the link diagram is \textit{classical}; otherwise, it is
\textit{virtual}. Equivalence classes of oriented link diagrams under
the Reidemeister moves correspond to ambient isotopy classes of
embedded disjoint oriented circles in $S^3$ if the link diagrams are
classical and stable equivalence classes of embedded disjoint oriented
circles in thickened surfaces $S\times [0,1]$ if the diagrams are
virtual; see \cite{SK, FRS:racks, Kup}.

\begin{figure}[!ht] 
$$\includegraphics[height=1in]{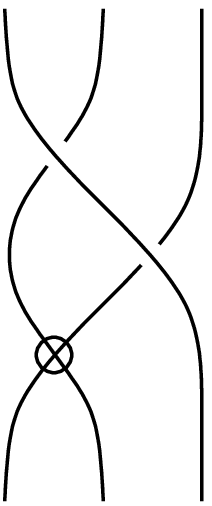}\quad\raisebox{0.5in}{$\iff$}\quad\includegraphics[height=1in]{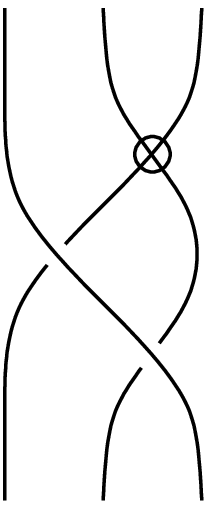}$$
\caption{The ``forbidden move'' $F_h$}\label{Fh}
\end{figure}

If the graph is non-planar, any additional crossings which must be
introduced in order to draw the graph in the plane are \textit{virtual
crossings}, self-intersections which are circled to distinguish them
from real crossings.  Any arc containing only virtual crossings may be
replaced by any other arc with the same endpoints and only virtual
crossings in a \textit{detour move} to obtain an equivalent virtual
link diagram. The equivalence relation on link diagrams generated by
the three Reidemeister moves and the detour move is called
\textit{virtual isotopy}. If, in addition, a strand with two classical
overcrossings is permitted to move past a virtual crossing, figure
\ref{Fh}, we have the \textit{forbidden move} $F_h$, and the
equivalence relation generated by virtual isotopy and $F_h$ is known
as \textit{welded isotopy} or \textit{weak virtual isotopy}; see
\cite{FRR, K}.  Welded links (ie welded isotopy classes of virtual
link diagrams) are important because of their close connection with
the braid-permutation group \cite{FRR}.

An oriented link diagram is a link diagram in which the edges are oriented 
in a coherent manner; specifically, of the two edges comprising the 
overstrand according to the crossing information, one is oriented toward the 
vertex and the other away from the vertex, and similarly for the undercrossing 
strand. These arcs are then unions of edges meeting at overcrossings.

Crossings are given signs according to their local writhe number. A crossing 
is \textit{positive} if the orientation on the plane determined by the 
orientations on the overarc followed by the underarc(s) gives the 
standard right-hand orientation on the plane; otherwise, the crossing is 
\textit{negative}.

In \cite{J} and \cite{FR}, combinatorial rules are given for associating a 
quandle to an oriented link diagram, with a generator for every arc and a 
relation at each crossing. In \cite{J}, the link diagrams are unoriented, 
with the relation determined by the blackboard framing of the link diagram, 
while in \cite{FR} the orientation of the over-crossing strand (but the not 
the under-crossing strand) determines the relation. Specifically, if we look 
in the positive direction of the overstrand $y$, we obtain the relation 
$x\tr y \sim z$ where $x$ is the undercrossing edge on the right and $z$ is
the undercrossing edge on the left. The quandle presented in this way is the
\textit{knot quandle} of the diagram; the fact that the knot quandle is an 
invariant of oriented knots and links is easily checked by comparing the 
presentations determined by the diagram before and after
each of the Reidemeister moves, keeping the quandle axioms in mind.

\begin{figure}[!ht] \small
\psfrag{x}{$x$}
\psfrag{y}{$y$}
\psfrag{z}{$z$}
\psfraga <-10pt,10pt> {+}{$+$}
\psfraga <-10pt,10pt> {-}{$-$}
\[\begin{array}{cc} 
\includegraphics[height=1in]{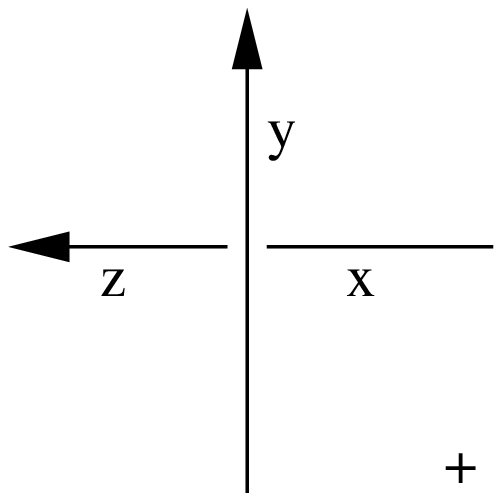}\qquad & \qquad \includegraphics [height=1in]{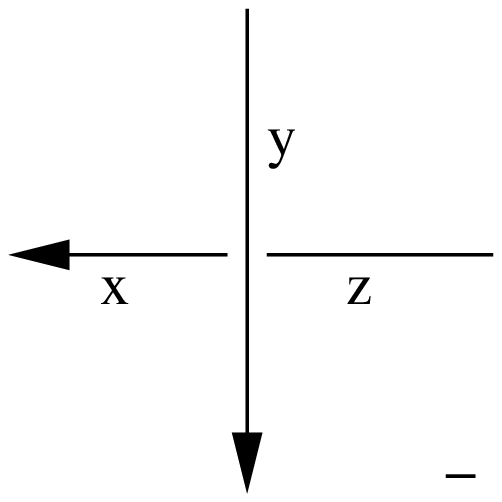} \\
 z\sim x\tr_+ y\qquad  &\qquad  x\sim z \tl_- y \\
 z\tl_+ y\sim x \qquad &\qquad  x\tr_- y \sim z
\end{array} \]
\caption{Crossing signs and quandle operations}\label{sign}
\end{figure}

The fact that $x\tr y\sim z$ is equivalent to $x\sim z\tl y$ shows that this 
rule determines not one unique relation but a pair of equivalent relations. 
This poses a difficulty when we wish to attempt to reconstruct a link diagram 
from its presented quandle; a short relation $x\sim z\tr y$ determines 
everything about a crossing except the orientation of the undercrossing 
strand. We can avoid this difficulty by including the crossing sign as
a subscript with each quandle operation.

It is then easy to see that this subscript convention yields information 
about which Tietze moves translate to Reidemeister moves and which do not. 
For instance, in every type II move, the two crossings have opposite signs,
whereas if we rewrite a short form relation with the other quandle operation,
both operations must have the same sign. This extra information is invisible
to the quandle structure and thus defines an enhanced quandle presentation,
called a \textit{signed quandle presentation}.

Not every combination of sign choice for each short form relation corresponds 
to a virtual link diagram; we must have exactly one inbound (that is, oriented 
toward the vertex) and one outbound (oriented away from the vertex) 
undercrossing arc at each crossing in order to have a virtual knot or link 
diagram. With our reconstruction rule, if the quandle operation is $\tr_+$ 
the input operand is inbound and the output operand is outbound. Switching 
either the sign or the triangle switches which operand is inbound and which 
is outbound. Simultaneously switching all of the signs gives us the quandle 
of the reflection of the virtual link, while switching the signs of a 
proper subset of relations in a presentation corresponding to a diagram 
yields an \textit{incoherent} signed quandle presentation -- one which does 
not correspond to a virtual link diagram. 

Indeed, since a choice of sign for one crossing determines which operand is 
inbound and which is outbound, such a choice also determines which operands 
are inbound and outbound in the other relations in which these generators 
appear are operands, and thus determines a sign for the crossings meeting 
the other ends of the undercrossing arcs. Likewise, each of these crossing 
signs determines the signs of the other crossings containing these operands,
and so on. In all, there are $2^N$ possible coherent quandle presentations 
for a given set of short relations, where $N$ is the number of components 
with at least one undercrossing.

Using the rule that distinct generators are assigned to every arc in an 
oriented knot diagram, we notice that every arc is either a simple closed 
curve or has endpoints at crossings; hence in a knot quandle presentation, 
every generator is either an operator-only generator or appears as an 
inbound operand exactly once and as an outbound operand exactly once in 
the set of relations.  Moreover, every relation determined by a 
crossing is in short form, and it is not hard to check that the remaining 
conditions in definition \ref{sform} are fulfilled.

Thus, the quandle presentation determined by a virtual link diagram is a 
short form quandle presentation. The converse is not true in general, since 
generators can appear more than twice as operands in a general short-form 
quandle presentation -- for example, 
\[\langle \ x,y,z, w \ | \ x\tr y \sim z, y\tr z \sim x, z\tr x \sim y, 
x\tr w \sim z \ \rangle.\]

\begin{definition}
A coherent short form signed quandle presentation is \textit{knotlike} if 
every generator either appears in exactly one relation as an inbound operand 
and in exactly one relation as an outbound operand, or appears only as an 
operator. A quandle is \textit{knotlike} if it has a knotlike presentation.
\end{definition}

\begin{figure}[!ht]\small
\psfrag{a}{$a$}
\psfraga <0pt, -2pt> {b}{$b$}
\[ \includegraphics{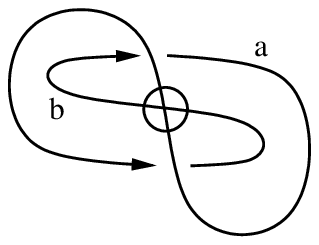} \quad \quad \includegraphics{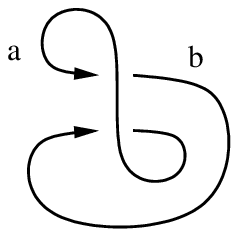} \]

\[Q = \langle a,b \ | \  a\tl_- a\sim b,  b\tl_- a \sim a \rangle\]

\caption{Example of non-isotopic virtual knots constructed from the same 
knotlike quandle presentation} \label{example1} 
\end{figure}

\begin{figure}[!ht]\small
\psfrag{a}{$a$}
\psfrag{b}{$b$}
\psfraga <-2pt, 0pt> {c}{$c$}
\psfrag{d}{$d$}
\[ \includegraphics{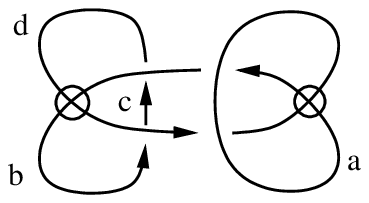} \quad \quad \includegraphics{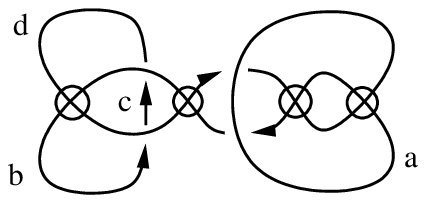} \]

\[Q = \langle a,b,c,d \ | \  b\tr_-a\sim a, \ b\tr_+ d\sim c, \ d\tr_-b \sim c,
\ d\tr_+a\sim a \rangle\]

\caption{Another example of non-isotopic virtual knots constructed from the 
same knotlike quandle presentation} \label{example2} 
\end{figure}

It is clear that every virtual link diagram defines a knotlike quandle
presentation.  Furthermore is it easy to see that every knotlike
quandle presentation determines a virtual link diagram; however a
little investigation shows that in most cases there are several
distinct virtual knot diagrams which all define the same signed
quandle presentation, eg, the inequivalent virtual knots pictured in
figures \ref{example1} and \ref{example2}. In light of this, it is
natural to ask to what degree a virtual link is determined by a
knotlike quandle presentation.

We can begin constructing a virtual link diagram from a knotlike quandle 
presentation by interpreting each signed short form relation as specifying a 
crossing, according to the reconstruction rule. The condition that every 
generator either appears exactly one relation as an inbound operand and 
exactly one relation as an outbound operand or only in the operator position
means that every arc either has a well-defined initial and terminal point
or is a simple closed curve, and we can complete the diagram by joining 
initial and terminal points with the same label with an arc including all 
the overcrossings of the same label, making virtual crossings as necessary.

The virtual knot diagram so constructed is not unique if any generator 
appears more than once as an operator. Indeed, for every such generator,
we must choose an order in which to go through the overcrossings along the 
arc, and different choices can result in non-isotopic virtual knot diagrams, 
all of which necessarily have the same quandle. Figure \ref{example1} depicts 
the two simplest inequivalent virtual knot diagrams constructed from the 
same quandle; the one on the left is Kauffman's virtualized trefoil, which 
is known to be non-trivial \cite{K}, while the diagram on the right is an 
unknot. 

Figure \ref{example2} shows another pair of known distinct virtual knots
constructed from the same knotlike quandle presentation, differing only in 
the order information; the virtual knot on the left in one of the Kishino 
virtual knots, known to be non-trivial, while the one on the right again is 
an unknot.

Figures \ref{example1} and \ref{example2} show that signed quandle
presentations alone are not sufficient to distinguish all oriented
virtual knots.  In order to completely specify a virtual knot diagram,
we must find a way of indicating an order for the overcrossings. In an
oriented knot diagram, each arc has a well-defined direction which we
may use to order the crossings; indeed, this is precisely what one
does when describing a virtual knot diagram by means of a Gauss
code. Moreover, since every crossing includes precisely one incoming
undercrossing strand, there is a bijection between the set of
generators and the set of crossings. Thus, the order of the inbound
undercrossings along a given arc in a virtual link diagram defines a
partial order on the set of generators, with two generators being
comparable if and only if the terminal points of their arcs lie along
the same overcrossing arc.  Call the sets of generators whose terminal
points lie along the the same arc \textit{arc-comparison classes}.

In addition to ordering crossings along arcs, a virtual knot diagram also 
defines a cyclic order in which the arcs are encountered while moving along 
a component. If the diagram is a virtual link, this cyclic ordering of arcs 
is a cyclic partial order on the arc-comparison classes of generators; if 
the diagram is a virtual knot, we have a cyclic order on these ordered sets.

\begin{definition} \label{partialorder}
Let $K$ be a virtual knot diagram and $Q=\langle X\ |\ R\rangle$ the 
corresponding knotlike quandle presentation. Say that $a<b$ if the terminal 
points of the arcs labeled $a$ and $b$ occur in crossings along the same 
arc $z$ and if traveling along $z$ from its initial point to its 
terminal point we encounter the terminal point of $a$ before the terminal 
point of $b$. 

Let $C_a$ be the arc-comparison class of $a$ under $<$. Then we have a cyclic
partial ordering $\prec$ on the set $\{C_x \ |\  x\in X\}$ defined by 
$C_a\prec C_b$ if the terminal point of $a$ meets the initial point of $b$
at a crossing. If the class $C_a$ consists of generators $x_1<x_2<\dots <x_n$,
we may write $C_a = \underline{x_1x_2\dots x_n}_a$. If $C_a$ is empty, we may
write $C_a = \underline{\qquad }_a$.

If $C_{a_1}\prec C_{a_2}\prec \dots \prec C_{a_n} \prec C_{a_1}$, we may simply
write $(C_{a_1}C_{a_2}\dots C_{a_n})$, and if $K$ has multiple components, 
we may specify the \textit{order information} of $K$ with a comma-separated 
list of these cyclically ordered partial order classes. 

Alternatively, we may specify the the order information implicitly by
simply listing the relations in the order specified by the order information,
with each relation's position determined by its inbound operand. The comparison
class of a generator $a$ then includes all relations with operator $a$.
\end{definition}

Since a generator $a$ is comparable to $b$ iff the terminal points of both
corresponding arcs lie along the same overcrossing arc, $<$ is a partial
order on the set of generators. Moreover, $<$ is clearly strict. If a 
generator $z$ is operator-only, its class $C_z$ will form its own class under 
$\prec$; in particular, the partial order $<$ on $C_z$ is actually a cyclic 
order on $C_z$ in this case, and only in this case -- the arc-comparison class
of any non-operator-only generator has a maximal element and minimal element, 
or is empty.

\begin{figure}[!ht]\small
\psfrag{a}{$a$}
\psfrag{b}{$b$}
\psfrag{c}{$c$}
\psfraga <0pt, -1pt> {d}{$d$} 
\psfrag{e}{$e$}
\psfraga <-2pt, 0pt> {f}{$f$}
\psfrag{g}{$g$}
\psfrag{h}{$h$}
\psfrag{j}{$j$}
\psfraga <-2pt, 0pt> {i}{$i$} 
\[\includegraphics{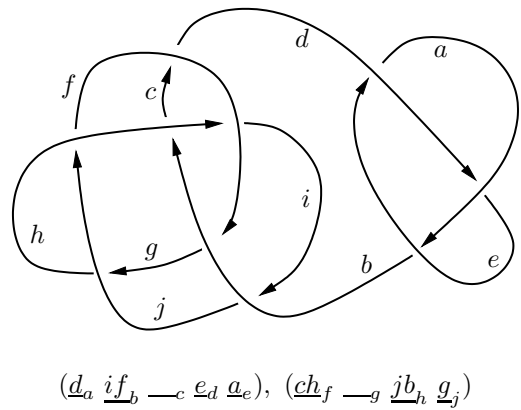}\]

\begin{center}
$(\underline{d}_a \ \underline{if}_b \ \underline{\quad }_c \ \underline{e}_d 
\ \underline{a}_e),  \
(\underline{ch}_f\  \underline{\quad }_g\  \underline{jb}_h\ \underline{g}_j)$
\end{center}
\caption{An example of order information from a link diagram}
\label{orderexample}
\end{figure}

From a virtual knot digram we can read off this order information and include
it with a knotlike quandle presentation to obtain a \textit{signed ordered 
knotlike quandle presentation} or a \textit{SOKQ presentation}. Conversely, 
given a signed ordered knotlike quandle presentation, we can construct a 
virtual link diagram which is uniquely determined up to virtual moves.

\section{Formal Reidemeister moves}\label{sec:frm}

In this section we note how Reidemeister moves on a virtual knot diagram
change an ordered quandle presentation. 

Reidemeister type I and II moves come in two varieties, crossing-introducing 
and crossing-removing. The crossing-introducing Reidemeister type I move 
breaks an arc into two and introduces a crossing, while the crossing-removing 
type I move removes a crossing and joins two arcs which previously met at an 
undercrossing into one.  The crossing-introducing type II move introduces
two crossings and breaks one arc into three, while the crossing-removing type 
II move removes two crossings and joins three arcs into one.

Let us use the convention that in a type I or II move, the original generator 
name stays with the terminal point of the arc. Then the crossing-introducing
type I move
\begin{list}{$\bullet$}{}
\item breaks the arc-comparison class $C_a=\underline{x_1\dots x_n}_a$ 
into two classes, either $C_{a'}=\underline{x_1\dots x_i a'}_{a'}$ and 
$C_a=\underline{x_{i+1}\dots x_n}_a$ or 
$C_{a'}=\underline{x_1\dots x_i}_{a'}$ and
$C_a=\underline{a'x_{i+1}\dots x_n}_a$,
\item introduces a new generator $a'$ and relation $a \sim a'\tr_+ a'$ or 
$a \sim a'\tl_- a'$ in the first case or $a'\sim a\tr_- a$ or 
$a'\sim a\tl_+ a$ in the second case, and 
\item replaces the operator $a$ with $a'$ in the relations with inbound 
operand $x_1,\cdots x_i$ and replaces the outbound operand $a$ with $a'$ 
in its relation.
\end{list}

Conversely, a crossing-removing type I move is only available in an 
ordered quandle presentation if the arc-comparison classes $C_a' \prec C_a'$
are adjacent in the cyclic order, with $C_{a'}$ ending in $a'$ if the relation 
is $a \sim a'\tr_+ a'$ or $a \sim a'\tl_- a'$, or $C_a$ starting with $a'$ if
the relation is $a'\sim a\tr_- a$ or $a'\sim a\tl_+ a$. In this case, we 
can delete the generator $a'$ and the relation with inbound operand $a'$, 
replace every instance of $a'$ with $a$, and join $C_{a'}$ with $C_a$ in 
order, minus the $a'$, to form the new $C_a$ and thus obtain the new signed
ordered quandle presentation. See figure 6 for an illustration.

Either of these operations on an ordered quandle presentation will be
called \textit{formal Reidemeister move~I\/}. Note that in terms of
Tietze moves, we've simply introduced a new generator $a'$ and
defining relation $a'\sim a$, replaced some instances of $a$ with $a'$
and made the presentation knotlike by replacing $a\sim a'$ with an
equivalent relation.  Absent the order information, we can replace
any, all, or none of the operator occurrences of $a$ with $a'$ to
obtain an equivalent quandle, while most such moves will not
correspond to Reidemeister moves or move sequences, though many may be
realizable as welded isotopy sequences.

\begin{figure}[!ht] \small
\psfraga <-3pt, 7pt> {a}{$a$}
\psfrag{b}{$b$}
\psfrag{c}{$c$}
\psfrag{y}{$y$}
\psfrag{z}{$z$}
\psfrag{x1}{$x_1$}
\psfrag{xi}{$x_i$}
\psfrag{xj}{$x_j$}
\psfrag{xn}{$x_n$}
\begin{center}
\begin{tabular}{ccc}
\includegraphics[height=2in]{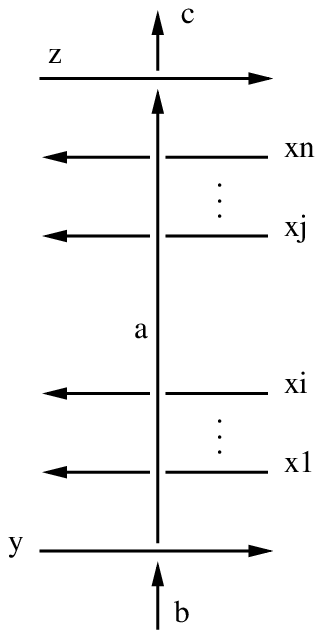} &  \hskip 0.3in \raisebox{1in}{$\iff$} &
\includegraphics[height=2in]{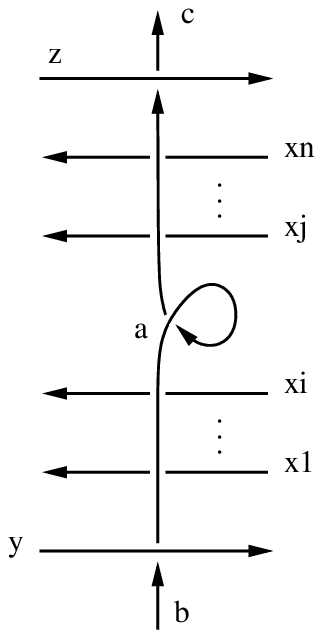} \\
$b\tr_+ y \sim a,\ a\tr_+ z\sim c$ & & $b\tr_+ y \sim a',\ a'\tr_+a' \sim a,\ 
a\tr_+ z\sim c$ \\
$\underline{x_1,\dots,x_n}_a$ & &  $\underline{x_1,\dots,x_ia'}_{a'}\ 
\underline{x_{i+1},\dots,x_n}_{a}$ \\
 & & \\
 & & operator $a$ replaced with \\
 & & $a'$ in relations involving \\
 & & inbound operand $x_1\dots x_i$
\end{tabular}
\end{center}
\caption{Formal Reidemeister move I example}\label{RI}
\end{figure} 

A crossing-introducing Reidemeister II move breaks an arc-comparison
class into three classes while introducing two generators, two relations,
and replacing some instances of one generator with a new one. Specifically,
\begin{list}{$\bullet$}{}
\item the arc-comparison class $C_a=\underline{x_1\dots x_n}_a$ is replaced 
by $C_{a''}C_{a'}C_{a} =\underline{x_1\dots{x_i}}_{a''}\underline{\quad}_{a'}
\underline{x_{i+1}\dots x_n}_a$ and the generators $a''a'$ or $a'a''$ are 
inserted somewhere in $C_z$,
\item new generators $a''$ and $a'$ are introduced along with relations 
$a'\sim a \tl_- z$ and $a''\sim a'\tr_+ z$ or
$a'\sim a \tr_- z$ and $a''\sim a'\tl_+ z$, and
\item the outbound operand $a$ is replaced in its relation with $a''$ and the 
operator $a$ is replaced with $a''$ in the relations with inbound operands 
$x_1,\dots, x_i$.
\end{list}

Conversely, a crossing-removing Reidemeister II move is available only when we 
have relations with opposite crossing signs 
$a'\sim a \tl_- z$ and $a''\sim a'\tr_+ z$ or
$a'\sim a \tr_- z$ and $a''\sim a'\tl_+ z$,
and order information including
\[ C_{a''}C_{a'}C_{a} =\underline{x_1\dots{x_i}}_{a''}\underline{\quad}_{a'}
\underline{x_{i+1}\dots x_n}_a, \]
and\[ 
C_z=\underline{\dots a'a'' \dots}_z \ \ \mathrm{or} \ \ 
C_z=\underline{\dots a''a' \dots}_z;\] 
in this case we may delete the generators $a'',a'$, their relations and 
replace any remaining instances of $a''$ with $a$.

\begin{figure}[!ht] \small
\psfrag{a}{$a$}
\psfrag{a1}{$a_1$}
\psfrag{a2}{$a_2$}
\psfrag{b}{$b$}
\psfrag{c}{$c$}
\psfrag{w}{$w$}
\psfrag{z}{$z$}
\psfrag{x1}{$x_1$}
\psfrag{xi}{$x_i$}
\psfrag{xj}{$x_j$}
\psfrag{xn}{$x_n$}
\psfrag{zj}{$z_j$}
\psfraga <-6pt,-2pt> {zk}{$z_{j+1}$}
\begin{center}
\begin{tabular}{ccc}
\includegraphics{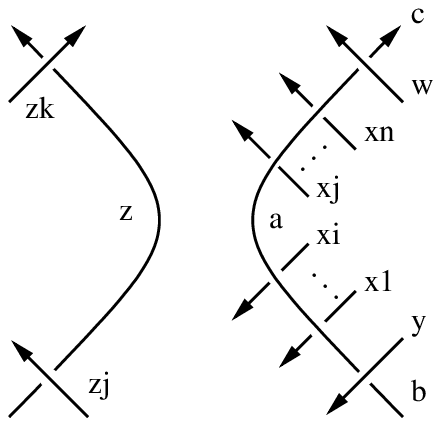} & \raisebox{1in}{$\iff$} &
\includegraphics{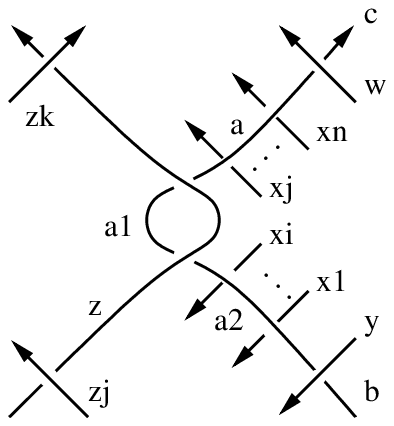} \\
$b\tr_+ y\sim a$
& & $b\tr_+ y\sim a''$,
$a''\tr_+ z \sim a'$, \\
& & $a' \tl_- z \sim a$ \\
$\underline{\dots z_jz_{j+1}\dots}_z,$ $\underline{x_1,\dots,x_n}_a$ & & 
$\underline{\dots z_ja''a'z_{j+1}\dots}_z,$ \\
& &  $\underline{x_1,\dots,x_i}_{a''}\ 
\underline{\quad}_{a'}\ \underline{x_{i+1},\dots,x_n}_{a}$ \\ \\
& & operator $a$ replaced with  \\
& & operator $a''$ in relations \\
& & with inbound operand $x_1,\dots, x_i$
\end{tabular}
\end{center}
\caption{Formal Reidemeister move II example}\label{RII}
\end{figure} 

In the Reidemeister type III move, one generator is replaced with another, 
and since this generator does not appear as an operator in any relation,
for simplicity we may use the same name for both the old and new generator. 
There are a number of cases, but in each case, we have a set of three
short relations which get replaced by another set of three short 
relations with one relation the same, one (the defining relation for
the generator which gets removed and re-added) changed, and the other
relation changed by a Tietze move involving a right-distribution.

The order information in a type III move changes by a cyclic permutation
of the input operands around the central triangle formed by the three strands
in the move. Since the cyclic order of the strands does not change in the 
move, the cyclic order of the arc-comparison classes also does not change. 
The arc-comparison class corresponding to the top strand contains 
two of the pictured generators, the arc-comparison class of one end of the 
middle strand has one pictured generator either as its maximum or minimum 
element, and the arc-comparison class of the other end of the middle strand 
includes no pictured generators. The arc-comparison classes of the bottom 
strand generators are unaffected by the move.

The type III move fixes the cyclic order of the arc-comparison classes, and it 
changes the partial order of the generators inside the classes in a nice way.
The two arc-comparison classes with operators on the middle strand are adjacent
in the cyclic order, with the pictured generator at the end of one class next 
to the other. The other two generators with pictured terminal points are 
adjacent somewhere in the arc-comparison class of the top strand. The move 
changes both the positions occupied by the pictured generators and which 
generators fill those positions. The position occupied by the pair remains in 
place, while the position at the extreme end of one of the middle-strand 
classes moves to the other extreme of the adjacent class. The generators 
filling these positions then undergo a cyclic permutation.

\begin{figure}[!ht]\small
\psfraga <1pt, 1pt> {u1}{$u_1$}
\psfrag{u}{$u$}
\psfrag{v}{$v$}
\psfrag{t}{$t$}
\psfrag{z}{$z$}
\psfrag{x}{$x$}
\psfrag{y}{$y$}
\psfrag{yn}{$y_n$}
\psfraga <0pt, 1pt> {zi}{$z_i$}
\psfrag{zj}{$z_{i+1}$}
\[
\begin{array}{ccc}
\includegraphics{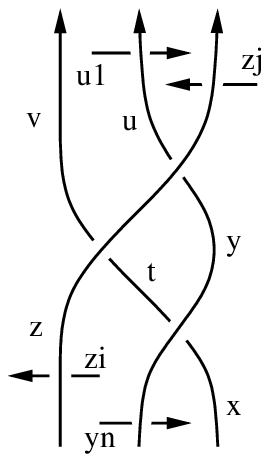} & \raisebox{1in}{$\iff$} &
\includegraphics{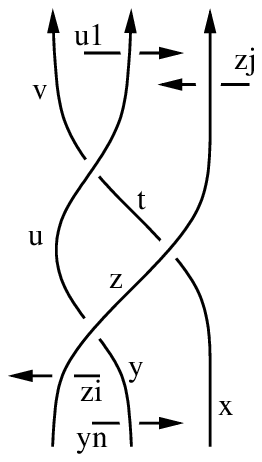} \\
y\tr_+ z \sim u & & y\tr_+ z \sim u  \\
t\tr_+ z \sim v & & t\tr_+ u \sim v \\
x\tr_+ y \sim t & & x\tr_+ z \sim t  \\
 & & \\
\underline{\dots z_i\ ty\ z_{i+1}\dots }_{z}, \underline{\dots y_nx}_{y} 
\underline{u_1 \dots }_{u} 
 & &
\underline{\dots z_i\ yx\ z_{i+1}\dots }_{z}, \underline{\dots y_n}_{y} 
\underline{tu_1\dots }_{u} 
\end{array}
\]
\caption{Formal Reidemeister move III example}
\end{figure} 

The fact that two oriented virtual knots are, by definition, virtually 
isotopic iff they are related by Reidemeister moves implies the following:

\begin{theorem}\label{main}
Two signed ordered knotlike quandle presentations present isotopic oriented
virtual knots if and only if they are related by a finite sequence of formal 
Reidemeister moves.
\end{theorem}

Note that in all three moves, the cyclic order of the arc-comparison classes
is changed only by insertions of new generators and deletions of old ones. In 
particular, the cyclic order of the arc-comparison classes is an invariant of 
oriented virtual knot type, in the following sense:

\begin{proposition} If 
$C_{x_1}\prec C_{x_2} \prec C_{x_3}$ in $Q=\langle X | R \rangle$, 
$Q'=\langle X' | R' \rangle$ and $f\co X\to X'$ 
generates an isomorphism of quandles which is realizable by formal 
Reidemeister moves, then the cyclic order on $Q'$ must include 
$C_{f(x_1)}\prec C_{f(x_2)} \prec C_{f(x_3)}$. 
\end{proposition}

While the cyclic order does not distinguish the two diagrams in figure
\ref{example1}, the full order information is different for these two. 
Specifically, the virtual knot diagram on the left has order information
$\underline{ab}_a\underline{\quad}_b$, while the diagram on the right
has order information $\underline{ba}_a\underline{\quad}_b$. Similarly, the
order information for the non-trivial Kishino virtual knot on the left in 
figure \ref{example2} differs from that of the unknot on the right only by 
switching the order of the generators $a$ and $d$ in the arc-comparison set 
$C_a$; the diagram on the left has order information $\underline{ad}_a 
\underline{c}_b \underline{\quad}_c \underline{b}_d$, while the unknot 
diagram on the right has
order information $\underline{da}_a\underline{c}_b\underline{\quad}_c 
\underline{b}_d$. 

\begin{figure}[!ht]\small 
\psfrag{a}{$a$}
\psfrag{b}{$b$}
\begin{center}
\begin{tabular}{cc}
\includegraphics{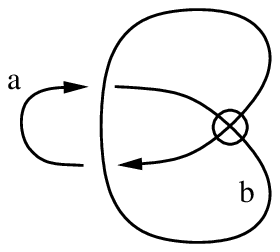} & \includegraphics{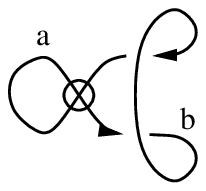} \\
$\langle a,b \ | \ b\tr_+ b \sim a, b\tr_-b \sim a, 
\underline{\quad}_a \underline{ba}_b \rangle$ & 
$\langle a,b \ | \ b\tr_+ b \sim a, b\tr_-b \sim a, 
\underline{\quad}_a \underline{ab}_b \rangle$
\end{tabular}
\end{center}
\caption{Two isotopic virtual knots whose SOKQ presentations differ only 
in the order information}\label{c}
\end{figure}

However, it not the case that SOKQ presentations which differ
only in the order information necessarily present non-isotopic virtual 
knots. For example, figure \ref{c} shows two isotopic virtual knots 
constructed from SOKQ presentations which differ only in the order 
information. Moreover, a connected sum of the two virtual knots in 
this example with any other virtual knot will yield additional examples
of possibly isotopic virtual knots with the same quandle but different order 
information. Indeed, compare figure \ref{example2}.
In particular, the order information itself is not an
invariant of virtual knot type; rather, it is the equivalence class
of order-information posets under the three Reidemeister moves
which is an invariant of virtual isotopy.

The fact that these examples are non-classical, however, 
leaves open the question of whether there exist isotopic classical
knots with SOKQ presentation differing only in the order information.

SOKQ presentations give a method for systematically listing all possible 
virtual knot diagrams with $n$ arcs. Namely, for a set of $n$ generators, 
choose ordered pairs of generators to be inbound and outbound operands.
Then for each such pair, select an operator, quandle operation and coherent 
crossing sign, determining a knotlike quandle presentation. For each resulting
signed knotlike quandle presentation, the operators divide the set of inbound 
generators into arc-comparison classes, and we then can list each possible 
ordering of the generators within the classes to complete an SOKQ 
presentation. For example, the two diagrams in figure \ref{example1} are the
only two possible diagrams, up to virtual moves, with the given knotlike 
quandle presentation. Note that the cyclic ordering of the classes is 
determined by the inbound/outbound pairing, though we can choose any 
ordering of the crossings within each arc class. 

We can now ask what other invariants of virtual isotopy might be determined 
by the order information in a signed ordered knotlike quandle presentation;
such an invariant is necessarily not determined by the quandle. The
primary new invariant contributed by SOKQ presentations is the equivalence
class of cyclically-ordered posets under the three formal Reidemeister
moves. It is not yet clear to this author how to compare two such 
poset-classes, though the problem deserves further study.

Other possibilities include invariants derived from the poset-class. 
For example, we have the following:

\begin{definition}
Let $K$ be an oriented virtual knot diagram. Label the arcs in $K$ with the 
numbers $1,2,\dots, n$ in the order in which they are encountered traveling 
around the 
knot in the direction of the orientation, so that the cyclic order of the 
arc comparison classes is $C_1\prec C_2 \prec \dots \prec C_n$. The endpoint 
ordering determines a permutation of $\{1,2,\dots, n\}$, specifically an 
$n$--cycle, which we may regard as an element of the infinite cyclic group 
$\Sigma_{\infty}$. We then define the \textit{order-permutation group} of the 
knot $K$ to be the subgroup of $\Sigma_{\infty}$ generated by these 
permutations for all diagrams of the oriented virtual knot $K$. 
\end{definition}

The full order information tells us which other permutations we must 
include in the generating set of the order-permutation group: a type I move 
inserts a number $i$ at some point along the arc $C_{i-1}$ and increments 
the  arc labels $i,i+1, \dots, n$ by 1. This may be 
accomplished by composing the permutation corresponding to the pre-move 
diagram with the transposition $(x(n+1))$ where the insertion is after the 
symbol $x$, then composing with the cycle $(i(i+1)\dots (n+1))$ to increment
the arcs. A type II move is similar, inserting two adjacent 
numbers $i,i+1$ anywhere in the cycle and incrementing the arc labels
$i,\dots, n$ by 2, done by composing the pre-move cycle with $(x(n+1)(n+2))$ 
or $(x(n+2)(n+1))$ then with $(i(i+1)\dots (n+1)(n+2))$. Finally, type III 
moves simply compose the pre-diagram cycle with a 3--cycle determined by the 
order information as described previously. While this group is infinitely 
generated, perhaps some more computable finite invariant may be derivable
from it.

\section{Welded links}\label{sec:weld}

We now turn to welded links.  Looking at the forbidden move (figure
\ref{Fh}) we observe that this move leaves invariant the cyclic order
information but changes the ordering within the arc-comparison
classes.  Thus we would expect that welded links (ie welded isotopy
classes of virtual link diagrams) should correspond to signed knotlike
quandle presentations (ie without the order information).  We shall
see that this is indeed the case.  Let us call a signed knotlike
quandle presentation an \textit{SKQ presentation\/} for short.

We can define formal Reidemeister moves on SKQ presentations in an
analogous way to those for SOKQ presentations.  Indeed the definitions
are given by simply omitting the order information from the preceding
ones.  We say that SKQ presentations are \textit{equivalent\/} if they
are related by a sequence of formal Reidemeister moves.  A virtual
link diagram gives an SKQ presentation as we saw earlier and
furthermore Reidemeister moves, the detour move and the forbidden move
$F_h$ on a diagram all correspond to formal Reidemeister moves.  Thus
we have a natural map $\Phi$ from oriented welded links to equivalence
classes of SKQ presentations.

\begin{theorem}\label{welded}
The natural map $\Phi$ just described is a bijection. 
\end{theorem}

\begin{proof}
We can turn an SKQ presentation into an SOKQ presentation by choosing
an arbitrary order for the arc-comparison classes, noting that the
cyclic order of the classes is determined by the SKQ presentation, and
hence by theorem \ref{main} find a corresponding virtual link diagram.
Furthermore a change in the chosen order can be realised as a product
of adjacent interchanges each of which can be realised by an $F_h$
move.  Thus we have a map from SKQ presentations to welded links.  Now
a formal Reidemeister move on an SKQ presentation corresponds to a
formal Reidemeister move on corresponding SOKQ presentations and hence
we have a map $\Psi$ from equivalence classes of SKQ presentations to
welded links.  It can readily be seen that $\Phi$ and $\Psi$ are
inverse bijections.\end{proof}

\section{The framed case}\label{sec:framed}

We finish by observing that there is an analogous treatment for framed
virtual links and framed welded links.  Framed virtual links are
equivalence classes under Reidemeister moves II and III of 4--valent
graphs with vertices interpreted as crossings, or equivalently of
framed virtual link diagrams under Reidemeister moves II and III and
the detour move.  For framed welded links the forbidden move $F_h$ is
also allowed.  It is worth commenting that for framed classical links
it is necessary to include a double Reidemeister I move \cite[page
370]{FR}.  With the detour move available is it an easy exercise to
check that this double move is unnecessary.

The corresponding algebraic object for framed links is the {\it rack\/}
which has exactly the same definition as the quandle but with axiom
(qi) omitted.  We can define signed knotlike rack presentations and
signed ordered knotlike rack presentations (SKR and SOKR presentations
respectively) by copying the quandle definitions without change.  Then
completely analogous arguments show:

\begin{theorem}The are bijections between the sets of isotopy classes of 
framed virtual links (respectively framed welded links) and the sets
of equivalence classes of SOKR presentations (respectively SKR
presentations) under the\break equivalence generated by formal Reidemeister
moves II and III.\end{theorem}

\subsection*{Acknowledgments}

The author would like to thank Xiao-Song Lin and Jim Hoste, whose
conversations the author found invaluable during the preparation of
this paper, as well as Louis Kauffman and Scott Carter, whose comments
at a conference started the author thinking about knot reconstruction
from quandles. The author would also like to thank the referee and the
editor, whose comments and corrections improved the paper
considerably.

\Addresses\recd

\end{document}

%% file: agtout.tex
\def\ifplaintex{\expandafter\ifx\csname documentclass\endcsname\relax}

\def\gtp{{\mathsurround=0pt\it $\cal G\mskip-2mu$eometry \&\ 
$\cal T\!\!$opology $\cal P\!$ublications}}  

\def\recd{{\small Received:\qua\receiveddate\ifx\reviseddate\relax
\else\qquad Revised:\qua\reviseddate\fi\par}} 

\def\lognumber#1{\def\thelognumber{#1}}
\def\volumenumber#1{\def\thevolumenumber{#1}}
\def\volumeyear#1{\def\thevolumeyear{#1}}
\def\papernumber#1{\def\thepapernumber{#1}}
\def\pagenumbers#1#2{\def\startpage{#1}\def\finishpage{#2}}
\def\published#1{\def\publishdate{#1}}

\def\received#1{\def\receiveddate{#1}}
\def\revised#1{\def\reviseddate{#1}}
\def\accepted#1{\def\accepteddate{#1}}

\let\\\par\let\thelognumber\relax\let\thevolumenumber\relax
\let\thepapernumber\relax\let\thevolumeyear\relax\let\startpage\relax
\let\finishpage\relax\let\publishdate\relax\let\receiveddate\relax
\let\reviseddate\relax\let\accepteddate\relax\let\theasciititle\relax
\let\theasciiauthors\relax
\let\theasciiabstract\relax

\let\theasciiemail\relax

\ifplaintex
\font\logobig=cmssbx10 scaled 3836
\font\logomed=cmssbx10 scaled 2557
\else
\font\logobig=cmssbx10 scaled 4200
\font\logomed=cmssbx10 scaled 2800
\fi

\long\def\makeagttitle{   
\count0=\startpage
\agt\hfill      
\hbox to 45truept{\vbox to 0pt{\vglue -13truept{\logomed A\kern -.37em{\logobig 
T}\kern -.38em G}\vss}\hss}
\break
{\small Volume \thevolumenumber\ (\thevolumeyear)
\startpage--\finishpage\nl
Published: \publishdate}

\vglue .25truein

{\parskip=0pt\leftskip 0pt plus
1fil\def\\{\par\smallskip}{\Large\bf\thetitle}\par\medskip} \vglue
0.05truein

{\parskip=0pt\leftskip 0pt plus 1fil\def\\{\par}{\sc\theauthors}
\par\medskip}%
 
\vglue 0.03truein 

{\small\leftskip 25truept\rightskip 25truept{\bf Abstract}\stdspace\theabstract

{\bf AMS Classification}\stdspace\theprimaryclass
\ifx\thesecondaryclass\relax\else; \thesecondaryclass\fi\par
{\bf Keywords}\stdspace \thekeywords\par}\vglue 7truept

}   

\ifplaintex
\font\phead=cmsl9 scaled 950
\font\pnum=cmbx10 scaled 913
\font\pfoot=cmsl9 scaled 950
\headline{\vbox to 0pt{\vskip -4.5mm\line{\small\phead\ifnum
\count0=\startpage ISSN 1472-2739 (on-line) 1472-2747 (printed)
\hfill {\pnum\folio}\else\ifodd\count0\def\\{ }%
\ifx\theshorttitle\relax\thetitle\else\theshorttitle\fi\hfill{\pnum\folio}
\else\def\\{ and }{\pnum\folio}\hfill\ifx\theshortauthors\relax\theauthors
\else\theshortauthors\fi\fi\fi}\vss}}
\footline{\vbox to 0pt{\vglue 0mm\line{\small\pfoot\ifnum\count0=\startpage
\copyright\ \gtp\hfill\else
\agt, Volume \thevolumenumber\ (\thevolumeyear)\hfill\fi}\vss}}
\else
\font=cmsl9 scaled 1050
\font\lnum=cmbx10 
\font=cmsl9 scaled 1050
\makeatletter
\def\@oddhead{{\small\lhead\ifnum\count0=\startpage ISSN 1472-2739 
(on-line) 1472-2747 (printed)\hfill {\lnum\number\count0}\else\ifodd\count0
\def\\{ }\ifx\theshorttitle\relax \thetitle \else\theshorttitle\fi\hfill
{\lnum\number\count0}\else\def\\{ and }{\lnum\number\count0}
\hfill\ifx\theshortauthors\relax 
\theauthors\else\theshortauthors\fi\fi\fi}}\def\@evenhead{\@oddhead}
\def\@oddfoot{\small\lfoot\ifnum\count0=\startpage\copyright\ \gtp\hfill\else
\agt, Volume \thevolumenumber\ (\thevolumeyear)\hfill\fi}
\def\@evenfoot{\@oddfoot}
\makeatother
\fi
\let\maketitlepage\makeagttitle
\let\makeshorttitle\maketitlepage
\let\maketitle\maketitlepage

\newwrite\gtoutfile
\long\gdef\makeheadfile{  
{\def\\{, }\def\s{ }
\immediate\openout\gtoutfile head.xxx
\immediate\write\gtoutfile{Proxy-for: \ifx\theasciiauthors\relax
\theauthors\else\theasciiauthors\fi\s<\ifx\theasciiemail\relax\theemail\else\theasciiemail\fi>}
\immediate\write\gtoutfile{\noexpand\\}
\immediate\write\gtoutfile{Authors: \ifx\theasciiauthors\relax
\theauthors\else\theasciiauthors\fi}
{\def\\{ }\immediate\write\gtoutfile{Title: \ifx\theasciititle\relax
\thetitle\else\theasciititle\fi}}
\immediate\write\gtoutfile{Subj-class: GT or SG, GR etc}
\immediate\write\gtoutfile{MSC-class: \theprimaryclass\ifx\thesecondaryclass\relax\else, \thesecondaryclass\fi}
\immediate\write\gtoutfile{Journal-ref: Algebr. Geom. Topol. \thevolumenumber\s
(\thevolumeyear) \startpage-\finishpage}
\immediate\write\gtoutfile{Comments: Published by Algebraic and
Geometric Topology at}
\immediate\write\gtoutfile{\s\s\s  http://www.maths.warwick.ac.uk/agt/AGTVol\thevolumenumber/agt-\thevolumenumber-\thepapernumber.abs.html}
\immediate\write\gtoutfile{\noexpand\\}
\immediate\write\gtoutfile{}
\ifx\theasciiabstract\relax
\immediate\write\gtoutfile{\theabstract}\else
\immediate\write\gtoutfile{\theasciiabstract}\fi
\immediate\write\gtoutfile{}
\immediate\write\gtoutfile{\noexpand\\}
\immediate\write\gtoutfile{}
\immediate\closeout\gtoutfile}}  

\def\maketitlepage{\makeagttitle\makeheadfile}
\let\makeshorttitle\maketitlepage
\let\maketitle\maketitlepage